\documentclass{amsart}

\usepackage{amsthm}
\usepackage[leqno]{amsmath}
\usepackage{latexsym,amsfonts,amssymb}
\usepackage{hyperref}

\newcommand{\numberseries}{\mdseries}   

\newlength{\thmtopspace}                
\newlength{\thmbotspace}                
\newlength{\thmheadspace}               
\newlength{\thmindent}                  

\renewcommand{\subparagraph}{\vspace*{\thmbotspace}}

\newtheoremstyle{bfupright head,slanted body}
                {\thmtopspace}{\thmbotspace}
                {\slshape}{\thmindent}{\bfseries}{.}{\thmheadspace}
                {{\numberseries \thmnumber{(#2) }}\thmnote{#3}}

\newtheoremstyle{bfupright head,upright body}
                {\thmtopspace}{\thmbotspace}
                {\upshape}{\thmindent}{\bfseries}{.}{\thmheadspace}
                {{\numberseries \thmnumber{(#2) }}\thmnote{#3}}

\newtheoremstyle{bfit head,upright body}
                {\thmtopspace}{\thmbotspace}
                {\upshape}{\thmindent}{\upshape}{.}{\thmheadspace}
                {{\numberseries\thmnumber{(#2) }}
                {\bfseries\itshape\thmnote{\negthickspace#3}}}

\newtheoremstyle{it head,upright body}
                {\thmtopspace}{\thmbotspace}
                {\upshape}{\thmindent}{\upshape}{.}{\thmheadspace}
                {{\numberseries\thmnumber{(#2) }}
                {\itshape\thmnote{\negthickspace#3}}}

\newtheoremstyle{fixed bf head,slanted body}
                {\thmtopspace}{\thmbotspace}{\slshape}
                {\thmindent}{\bfseries}{.}{\thmheadspace}
                {{\numberseries \thmnumber{(#2) }}\thmname{#1}\thmnote{ (#3)}}

\newtheoremstyle{fixed bf head,upright body}
                {\thmtopspace}{\thmbotspace}{\upshape}
                {\thmindent}{\bfseries}{.}{\thmheadspace}
                {{\numberseries \thmnumber{(#2) }}\thmname{#1}\thmnote{ (#3)}}

\newtheoremstyle{numbered paragraph}
                {\thmtopspace}{\thmbotspace}{\upshape}
                {\thmindent}{\upshape}{}{0pt}
                {{\numberseries \thmnumber{(#2) }}}

\newtheoremstyle{unnumbered paragraph}
                {\thmtopspace}{\thmbotspace}{\upshape}
                {\parindent}{\upshape}{}{0pt}

\theoremstyle{bfupright head,slanted body}
\newtheorem{res}{}[section]             \newtheorem*{res*}{}

\theoremstyle{bfupright head,upright body}
               \newtheorem*{bfhpg*}{}

\theoremstyle{fixed bf head,slanted body}
\newtheorem{thm}[res]{Theorem}          \newtheorem*{thm*}{Theorem}
\newtheorem{prp}[res]{Proposition}      \newtheorem*{prp*}{Proposition}
\newtheorem{cor}[res]{Corollary}        \newtheorem*{cor*}{Corollary}
\newtheorem{lem}[res]{Lemma}            \newtheorem*{lem*}{Lemma}

\theoremstyle{fixed bf head,upright body}
       \newtheorem*{dfn*}{Definition}
      \newtheorem*{obs*}{Observation}
\newtheorem{rmk}[res]{Remark}           \newtheorem*{rmk*}{Remark}
\newtheorem{exa}[res]{Example}          \newtheorem*{exa*}{Example}

\theoremstyle{numbered paragraph}
\newtheorem{ipg}[res]{}

\theoremstyle{unnumbered paragraph}
\newtheorem{ipg*}{}

\newlength{\thmlistleft}        
\newlength{\thmlistright}       
\newlength{\thmlistpartopsep}   
\newlength{\thmlisttopsep}      
\newlength{\thmlistparsep}      
\newlength{\thmlistitemsep}     

\newcounter{prt}
\newenvironment{prt}{\begin{list}{\upshape (\alph{prt})}%
    {\usecounter{prt}%
      \setlength{\leftmargin}{\thmlistleft}%
      \setlength{\labelwidth}{\thmlistleft}%
      \setlength{\rightmargin}{\thmlistright}%
      \setlength{\partopsep}{\thmlistpartopsep}%
      \setlength{\topsep}{\thmlisttopsep}%
      \setlength{\parsep}{\thmlistparsep}%
      \setlength{\itemsep}{\thmlistitemsep}}}%
  {\end{list}}%

 \newcommand{\prtlbl}[1]{{\upshape(#1)}}

\newenvironment{prf*}[1][Proof]{%
  \begin{proof}[\bf #1]
    \setcounter{equation}{0}
    }
  {\end{proof}
}

\newcommand{\pgref}[1]{(\ref{#1})}

\newcommand{\thmref}[2][Theorem~]{#1\pgref{thm:#2}}
\newcommand{\corref}[2][Corollary~]{#1\pgref{cor:#2}}
\newcommand{\prpref}[2][Proposition~]{#1\pgref{prp:#2}}
\newcommand{\lemref}[2][Lemma~]{#1\pgref{lem:#2}}

\newcommand{\exaref}[2][Example~]{#1\pgref{exa:#2}}
\newcommand{\rmkref}[2][Remark~]{#1\pgref{rmk:#2}}
\newcommand{\secref}[2][Section~]{#1\ref{sec:#2}}

 \newcommand{\partpgref}[2]{(\ref{#1})\prtlbl{#2}}
 \newcommand{\partlemref}[3][Lemma~]{#1\partpgref{lem:#2}{#3}}

\newcommand{\thmcite}[2][?]{\cite[thm.~#1]{#2}}
\newcommand{\corcite}[2][?]{\cite[cor.~#1]{#2}}
\newcommand{\prpcite}[2][?]{\cite[prop.~#1]{#2}}
\newcommand{\lemcite}[2][?]{\cite[lem.~#1]{#2}}

\newcommand{\setof}[3][\;]{\{#1#2 \mid #3#1\}}
 \newcommand{\ZZ}{\mathbb{Z}}

\newcommand{\f}{\varphi}
\newcommand{\m}{\mathfrak{m}}
\newcommand{\n}{\mathfrak{n}}
\newcommand{\p}{\mathfrak{p}}
\newcommand{\q}{\mathfrak{q}}

\newcommand{\is}{\cong}
\newcommand{\qis}{\simeq}

\renewcommand{\le}{\leqslant}
\renewcommand{\ge}{\geqslant}

\newcommand{\onto}{\twoheadrightarrow}

 \newcommand{\xra}[2][]{\xrightarrow[#1]{\;#2\;}}
 \newcommand{\qra}{\xra{\;\qis\;}}

\newcommand{\Rhat}{\widehat{R}}
\newcommand{\Shat}{\widehat{S}}

\newcommand{\mapdef}[4][\rightarrow]{\nobreak{#2\colon #3 #1 #4}}
\renewcommand{\Im}[1]{\nobreak{\operatorname{Im}#1}}
\newcommand{\Ker}[1]{\nobreak{\operatorname{Ker}#1}}

\newcommand{\dif}[2][]{{\partial}^{#2}_{#1}}

\renewcommand{\H}[2][]{\operatorname{H}_{#1}(#2)}

\newcommand{\Shift}[2][]{\mathsf{\Sigma}^{#1}{#2}}

\newcommand{\SpecR}{\operatorname{Spec}R}
\newcommand{\dptR}{\operatorname{depth}R}
\newcommand{\dimR}{\operatorname{dim}R}

\newcommand{\E}[2][R]{\operatorname{E}_{#1}(#2)}

\renewcommand{\dim}[2][R]{\operatorname{dim}_{#1}#2}
\newcommand{\wdt}[2][R]{\operatorname{width}_{#1}#2}
\newcommand{\dpt}[2][R]{\operatorname{depth}_{#1}#2}

\newcommand{\id}[2][R]{\operatorname{id}_{#1}#2}
\newcommand{\pd}[2][R]{\operatorname{pd}_{#1}#2}

\newcommand{\Gfd}[2][R]{\operatorname{Gfd}_{#1}#2}
\newcommand{\Gid}[2][R]{\operatorname{Gid}_{#1}#2}
\newcommand{\Hom}[3][R]{\operatorname{Hom}_{#1}(#2,#3)}
\newcommand{\RHom}[3][R]{\operatorname{\mathbf{R}Hom}_{#1}(#2,#3)}
\newcommand{\Ext}[4][R]{\operatorname{Ext}_{#1}^{#2}(#3,#4)}

\newcommand{\tp}[3][R]{\nobreak{#2\otimes_{#1}#3}}
 \newcommand{\tpp}[3][R]{(\tp[#1]{#2}{#3})}
 \newcommand{\Ltp}[3][R]{\nobreak{#2\otimes_{#1}^{\mathbf{L}}#3}}
 \newcommand{\Ltpp}[3][R]{(\Ltp[#1]{#2}{#3})}
\newcommand{\Tor}[4][R]{\operatorname{Tor}^{#1}_{#2}(#3,#4)}

\newcommand{\Cat}[2]{{\mathsf{#2}}(#1)}
\newcommand{\D}[1][R]{\Cat{#1}{D}}
\newcommand{\A}[1][R]{\Cat{#1}{A}}
\newcommand{\B}[1][R]{\Cat{#1}{B}}

\makeatletter
\def\@nobreak@#1{\mathchoice%
  {\nobreakdef@\displaystyle\f@size{#1}}%
  {\nobreakdef@\nobreakstyle\tf@size{\firstchoice@false #1}}%
  {\nobreakdef@\nobreakstyle\sf@size{\firstchoice@false #1}}%
  {\nobreakdef@\nobreakstyle\ssf@size{\firstchoice@false #1}}%
  \check@mathfonts}%
\def\nobreakdef@#1#2#3{\hbox{{%
                    \everymath{#1}%
                    \let\f@size#2\selectfont%
                    #3}}}%
\makeatother

\newcommand{\x}{\mathbf{x}}

\newcommand{\supp}[2][R]{\operatorname{supp}_{#1}#2}

\numberwithin{equation}{res}

\begin{document}

\title[Transfer of Gorenstein dimensions along ring
homomorphisms]{Transfer of Gorenstein dimensions\\ along ring
  homomorphisms}

\author{Lars Winther Christensen}

\address{Lars Winther Christensen, Department of Mathematics and
  Statistics, Texas Tech University, Lubbock, TX 79409-1042, U.S.A.}

\email{lars.w.christensen@ttu.edu}

\urladdr{http://www.math.ttu.edu/{\tiny $\sim$}lchriste}

\author{Sean Sather-Wagstaff}

\address{Sean Sather-Wagstaff, Department of Mathematics, NDSU Dept \#
  2750, PO Box 6050, Fargo, ND 58108-6050, U.S.A.}

\email{Sean.Sather-Wagstaff@ndsu.edu}

\urladdr{http://math.ndsu.nodak.edu/faculty/ssatherw/}

\dedicatory{Dedicated with gratitude to Hans-Bj\o rn Foxby, our
  teacher and friend}

\date{10 August 2009}



\begin{abstract}
  A central problem in the theory of Gorenstein dimensions over
  commutative noetherian rings is to find resolution-free
  characterizations of the modules for which these invariants are
  finite. Over local rings, this problem was recently solved for the
  Gorenstein flat and the Gorenstein projective dimensions; here we
  give a solution for the Gorenstein injective dimension. Moreover, we
  establish two formulas for the Gorenstein injective dimension of
  modules in terms of the depth invariant; they extend formulas for
  the injective dimension due to Bass and Chouinard.
\end{abstract}

\maketitle


\section*{Introduction}

\noindent Gorenstein dimensions are homological invariants that are
useful for identifying modules and ring homomorphisms with good
homological properties. This paper is concerned with the Gorenstein
injective dimension and the Gorenstein flat dimension, denoted
$\Gid[]{}$ and $\Gfd[]{}$, respectively. These invariants are defined
in terms of resolutions by modules from certain classes, the
Gorenstein injective and the Gorenstein flat modules. See Section 1
for definitions.

Let $R$ be a commutative noetherian ring. It is frequently useful to
know that finiteness of the classical homological dimensions of an
$R$-module $M$ can be detected by vanishing of (co)homology. For the
injective dimension one has
\begin{equation*}
  \id{M}= \sup{\setof{j}{\text{$\Ext{j}{R/\p}{M}\neq 0$ for some $\p\in\SpecR$}}}.
\end{equation*}
One of the key problems in Gorenstein homological algebra has been to
find criteria for finiteness of Gorenstein dimensions that are
resolution-free.  See the survey \cite{CFH-} and the introduction in
\cite{CFH-06} for a further discussion of this issue. The problem was
partly solved by Christensen, Frankild, and Holm in~\cite{CFH-06}: If
$R$ has a dualizing complex and $M$ is an $R$-module, then
\begin{equation*}
  \text{$\Gid{M}$ is finite if an only if $M$ belongs to $\B$}
\end{equation*}
where $\B$ is the Bass class of $R$; the crucial point is that
verification of membership in $\B$ does not involve construction of a
Gorenstein injective resolution. Similarly $\Gfd{M}$ is finite if an
only if $M$ belongs to the Auslander class $\A$.

If $R$ is a homomorphic image of a Gorenstein ring, then it has a
dualizing complex. In particular \cite{CFH-06} solves the problem when
$R$ is local and complete or, more generally, essentially of finite
type over a complete local ring. However, non-trivial modules of
finite Gorenstein injective dimension or finite Gorenstein flat
dimension may exist over rings that are not homomorphic images of
Gorenstein rings; see \exaref{nodc}.  In \cite{MAEMTs07} Esmkhani and
Tousi show that when $R$ is local, but not necessarily a homomorphic
image of a Gorenstein ring, an $R$-module $M$ has finite Gorenstein
flat dimension if and only if the module $\tp{\Rhat}{M}$ is in the
Auslander class $\A[\Rhat]$ of the completion $\Rhat$.  This solves
the resolution-free characterization problem for the Gorenstein flat
dimension over local rings. In a separate paper~\cite{MAEMTs07a} the
same authors give a solution for the Gorenstein injective dimension of
cotorsion modules over local rings.

In this paper, we complete the solution for local rings with the
special case $S=\Rhat$ of the next result, wherein $\RHom{S}{M}$ and
$\Ltp{S}{M}$ are the right derived homomorphism complex and the left
derived tensor product complex.  More general statements are proved in
\pgref{thm:01} and~\pgref{thm:02}.

\begin{res*}[Theorem~A]
  Let $\mapdef{\f}{R}{S}$ be a local ring homomorphism such that $S$
  has a bounded resolution by flat $R$-modules when considered as an
  $R$-module via $\f$.  For every $R$-module $M$ there are
  inequalities
  \begin{equation*}
    \Gid{M} \ge \Gid[S]{\RHom{S}{M}} \quad\text{and}\quad
    \Gfd{M} \ge \Gfd[S]{\Ltpp{S}{M}}.
  \end{equation*}
  If $\f$ is flat, then equalities hold; in particular, the respective
  dimensions are simultaneously finite in this case.
\end{res*}

\noindent As noted above, Esmkhani and Tousi's~\cite{MAEMTs07a}
resolution-free characterization of finiteness of Gorenstein injective
dimension only applies to cotorsion modules.  The cotorsion hypothesis
is quite restrictive. Indeed, work of Frankild, Sather-Wagstaff, and
Wiegand~\cite{AJFSSW08, FSW-08} shows that a finitely generated
cotorsion $R$-module is complete.  For Gorenstein rings, the following
application of Theorem~A strengthens the main result of
\cite{AJFSSW08}; it only assumes that the Ext-modules are finitely
generated over $\Rhat$, not over $R$.

\begin{res*}[Theorem~B]
  Let $R$ be a Gorenstein local ring, and let $M$ be a finitely
  generated $R$-module. If the $\Rhat$-modules $\Ext{i}{\Rhat}{M}$ are
  finitely generated for $i=1,\dots,\dim{M}$, then the modules
  $\Ext{i}{\Rhat}{M}$ vanish for $i\ge 1$, and $M$ is complete.
\end{res*}

\noindent 
The hypotheses of this result are satisfied if $M$ is complete, e.g.,
if $M$ has finite length; cf.~\rmkref{no}. Theorem~B is a special case
of \thmref[]{FSWW}.

In \secref{bass} we consider formulas that express the Gorenstein
injective dimension of an $R$-module in terms of the depth
invariant. Our main result in this direction is Theorem~C below.  It
extends Chouinard's~\cite{LGC76} formula for injective dimension, and
it removes the assumption about existence of a dualizing complex from
\thmcite[6.8]{CFH-06}.

\begin{res*}[Theorem~C]
  For every $R$-module $M$ of finite Gorenstein injective dimension
  there is an equality
  \begin{displaymath}
    \Gid{M}=\sup\{\dpt[]{R_\p}-\wdt[R_\p]{M_\p} \,|\, \p \in \SpecR\}.
  \end{displaymath}
\end{res*}

\noindent
For certain modules this formula has already been established by
Khatami, Tousi, and Yassemi~\cite{KTY-,LKhSYs07}. Actually, we prove
Theorems~A and C for $R$-complexes, and the latter yields a Bass
formula for homologically finite $R$-complexes; see
\corref[]{Bass}. Such a formula was established for modules
in~\cite{KTY-}.

\section{Finiteness and descent of Gorenstein homological dimensions}
\label{sec:proofs}

\noindent Throughout this paper $R$ and $S$ are commutative noetherian
rings.  Complexes of $R$-modules, \emph{$R$-complexes} for short, are
indexed homologically: the $i$th differential of an $R$-complex $M$ is
written $\mapdef{\dif[i]{M}}{M_i}{M_{i-1}}$.  We proceed by recalling
the definitions of Gorenstein injective and Gorenstein flat modules
from \cite{EEnOJn95b, EJT-93}.

\begin{ipg}
  An $R$-module $J$ is said to be \emph{Gorenstein injective} if there
  is an exact complex $I$ of injective $R$-modules such that
  $J\cong\Ker{\dif[0]{I}}$ and the complex $\Hom{E}{I}$ is exact for
  every injective $R$-module $E$.

  An $R$-module $G$ is \emph{Gorenstein flat} if there is an exact
  complex $F$ of flat $R$-modules such that $G\cong\Im{\dif[0]{F}}$
  and $\tp{E}{F}$ is exact for every injective $R$-module~$E$.
\end{ipg}

The first step toward Theorem~A is to notice that the central
arguments in the works of Esmkhani and Tousi \cite{MAEMTs07,
  MAEMTs07a} apply to any faithfully flat ring homomorphism, not just
to the map $R \to \Rhat$; see \lemref[]{ET}. To this end the next fact
is key.

\begin{lem}
  \label{lem:injsum}
  Let $\mapdef{\f}{R}{S}$ be a faithfully flat ring homomorphism. If
  $E$ is an injective $R$-module, then it is a direct summand (as an
  $R$-module) of the injective $S$-module $\Hom{S}{E}$.
\end{lem}

\begin{prf*}
  Let $E$ be an injective $R$-module. It is well-known, and
  straightforward to show, that $\Hom{S}{E}$ is an injective
  $S$-module.  Because $\f$ is faithful, it is a pure monomorphism of
  $R$-modules, cf.~\thmcite[7.5]{Mat}. This implies that $S/R$ is a
  flat $R$-module, so $\Hom{S/R}{E}$ is injective. Now apply the exact
  functor $\Hom{-}{E}$ to $0 \to R \to S \to S/R \to 0$ to obtain a
  split exact sequence of injective modules.
\end{prf*}

\begin{lem}
  \label{lem:ET}
  Let $\mapdef{\f}{R}{S}$ be a faithfully flat ring homomorphism.
  \begin{prt}
  \item \label{item:ETb} Assume that $\dim[]{S}$ is finite. An
    $R$-module $M$ is Gorenstein injective if and only if\,
    $\Hom{S}{M}$ is a Gorenstein injective $S$-module and
    $\Ext{i}{F}{M}=0$ for every flat $R$-module $F$ and all $i \ge
    1$.\footnote[3]{ The vanishing of $\Ext{i}{F}{M}$ for every flat
      $R$-module $F$ and for all $i\geq 1$ means exactly that $M$ is
      \emph{cotorsion}. It is straightforward to show that this is
      equivalent to the standard definition of cotorsion which only
      requires $\Ext{1}{F}{M}=0$ for every flat $R$-module $F$.}

  \item \label{item:ETa} An $R$-module $M$ is Gorenstein flat if and
    only if $\tp{S}{M}$ is a Gorenstein flat $S$-module and
    $\Tor{i}{E}{M}=0$ for every injective $R$-module $E$ and all $i
    \ge 1$.
  \end{prt}
\end{lem}

\begin{prf*}
  Argue as in the proofs of \thmcite[2.5]{MAEMTs07} and
  \thmcite[2.5]{MAEMTs07a}, but use the in\-jective $S$-module
  $\Hom{S}{E}$ from \lemref[]{injsum} in place of the double Matlis
  dual $E^{\vee\vee}$.
\end{prf*}

For the proofs that follow, we need some terminology.

\begin{ipg}
  Let $M$ be an $R$-complex; it is said to be \emph{bounded above} if
  $M_i =0$ for $i \gg 0$, \emph{bounded~below} if $M_i =0$ for $i \ll
  0$, and \emph{bounded} if $M_i = 0$ for $|i| \gg 0$.  If the
  homology complex $\H{M}$ is bounded, then $M$ is called
  \emph{homologically bounded}. If $\H{M}$ is finitely generated, then
  $M$ is said to be \emph{homologically finite}. The notations
  $\inf{M}$ and $\sup{M}$ stand for the infimum and supremum of the
  set \mbox{$\{i\in\ZZ\mid\H[i]{M}\ne 0\}$}, with the convention that
  $\inf{M} = \infty$ and $\sup{M} = -\infty$ if $\H{M}=0$.

  From this point, we work in the derived categories $\D$ and $\D[S]$;
  see e.g.\ \cite{GelMan}.  Given two $R$-complexes $M$ and $N$, their
  left derived tensor product complex and right derived homomorphism
  complex are denoted $\Ltp{M}{N}$ and $\RHom{M}{N}$. The symbol
  `$\qis$' is used to identify isomorphisms in derived categories.
\end{ipg}

\begin{ipg}
  The \emph{Gorenstein injective dimension} of a homologically bounded
  $R$-complex $M$ is defined as follows
  \begin{equation*}
    \Gid{M} = \inf \left\{ \, \sup \{ i \in \ZZ \, | \, J_{-i} \ne 0 \} \,
      \left|
        \begin{array}{l}
          \mbox{$J$ is a bounded above complex} \\
          \mbox{of Gorenstein injective modules} \\
          \mbox{and isomorphic to $M$ in $\D$}
        \end{array}
      \right.
    \right\}.
  \end{equation*}
  The \emph{Gorenstein flat dimension} is defined similarly in terms
  of bounded below complexes of Gorenstein flat modules; see
  \cite[(5.2.3)]{lnm}.

  When $R$ has a dualizing complex $D$, Avramov and Foxby
  \cite{LLAHBF97} define two full subcategories $\A$ and $\B$ of $\D$.
  The objects in the \emph{Auslander class} $\A$ are the homologically
  bounded $R$-complexes $M$ such that $\Ltp{D}{M}$ is homologically
  bounded and the natural morphism \mbox{$M\to\RHom{D}{\Ltp{D}{M}}$}
  is an isomorphism in $\D$.  The objects in the \emph{Bass class}
  $\B$ are the homologically bounded $R$-complexes $M$ such that
  $\RHom{D}{M}$ is homologically bounded and the natural morphism
  \mbox{$\Ltp{D}{\RHom{D}{M}} \to M$} is an isomorphism in $\D$.
  
  For a homologically bounded $R$-complex $M$, the main results in
  \cite{CFH-06} state
  \begin{align}
    \label{A}
    \Gfd{M} &\text{ is finite if and only if $M$ belongs to $\A$;
      and}\\
    \label{B}
    \Gid{M} &\text{ is finite if and only if $M$ belongs to $\B$.}
  \end{align}
\end{ipg}

Before proving Theorem~A, we recall an elementary construction of
rings that admit non-trivial modules of finite Gorenstein dimensions.

\begin{exa}
  \label{exa:nodc}
  Let $Q$ be a commutative noetherian ring and consider the ring of
  dual numbers $R=Q[X]/(X^2)$.  It is routine to show that the cyclic
  $R$-module $R/(X)$ is Gorenstein flat and not flat. Hence, for every
  faithfully injective $R$-module~$E$ the module $\Hom{R/(X)}{E}$ is
  Gorenstein injective and not injective; see \thmcite[(6.4.2)]{lnm}.
  Furthermore, if $Q$ is not a homomorphic image of a Gorenstein ring,
  then neither is $R$.
\end{exa}

The next result contains half of Theorem~A from the
introduction. Recall that a ring homomorphism $\mapdef{\f}{R}{S}$ has
\emph{finite flat dimension} when $S$, considered as an $R$-module via
$\f$, has a bounded resolution by flat $R$-modules.

\begin{thm}
  \label{thm:01}
  Let $\mapdef{\f}{R}{S}$ be a ring homomorphism of finite flat
  dimension, and assume that $\dim[]{R}$ is finite.  For every
  homologically bounded $R$-complex $M$ there is an inequality
  \begin{align*}
    \Gid{M} &\ge \Gid[S]{\RHom{S}{M}}.
  \end{align*}
  If $\f$ is faithfully flat and $\dim[]{S}$ is finite, then equality
  holds; in particular, the dimensions are simultaneously finite in
  this case.
\end{thm}

\begin{prf*}
  Assume that $M$ has finite Gorenstein injective dimension, and fix a
  bounded complex $J$ of Gorenstein injective $R$-modules such that
  there is an isomorphism $M \qis J$ in $\D$. As an $R$-module, $S$
  has projective dimension at most $\dim[]{R}<\infty$; see
  \cite[II.~thm.~(3.2.6)]{LGrMRn71} and
  \cite[prop.~6]{CUJ70}. Therefore, by \corcite[2.12]{CFH-06} there is
  an isomorphism $\RHom{S}{M} \qis \Hom{S}{J}$ in $\D[S]$, and the
  right-hand complex is a bounded one of Gorenstein injective
  $S$-modules; see \cite[Ascent table II (h)]{LWCHHl09a}. In
  particular, there is an inequality $\Gid{M} \ge
  \Gid[S]{\RHom{S}{M}}$.

  Assume now that $\f$ is faithfully flat and that $d:=\dim[]{S}$ and
  $\Gid[S]{\RHom{S}{M}}$ are finite. Recall the inequalities $\pd{S}
  \le \dimR \le d$. Consider a resolution $M \xrightarrow{\simeq} I$
  by injective $R$-modules. The complex $\Hom{S}{I} \qis \RHom{S}{M}$
  is one of injective $S$-modules, and one has $\H[i]{\Hom{S}{I}}=0$
  for all $i<\inf{M} - d$, as $\pd{S}$ is at most $d$. Left-exactness
  of the functor $\Hom{S}{-}$ yields an isomorphism
  \begin{equation*}
    \Ker{\dif[n-1]{\Hom{S}{I}}} \is \Hom{S}{\Ker{\dif[n-1]{I}}}
  \end{equation*}
  for each $n$. It follows that $\RHom{S}{M}$ is isomorphic in $\D[S]$
  to the complex
  \begin{equation*}
    0 \to \Hom{S}{I_0} \to \dots \to
    \Hom{S}{I_{n}} \to \Hom{S}{\Ker{\dif[n-1]{I}}} \to 0
  \end{equation*}
  for $n<\inf{M}-d$. Set $K = \Ker{\dif[\inf{M}-2d-1]{I}}$. Since the
  $S$-complex $\RHom{S}{M}$ has finite Gorenstein injective dimension,
  the $S$-module $\Hom{S}{K}$ is Gorenstein injective; see
  \thmcite[3.3]{CFH-06}. To show that $\Gid{M}$ is finite, we
  use \partlemref{ET}{a} to prove that $K$ is Gorenstein injective
  over $R$: For every flat $R$-module $F$, one has $\pd{F}\le d$, and
  for every $i\le 1$ dimension shifting yields
  \begin{equation*}
    \Ext{i}{F}{K} \is \Ext{i+d}{F}{\Ker{\dif[\inf{M}-d-1]{I}}} =0.
  \end{equation*}
  To prove the equality of Gorenstein injective dimensions, choose an
  injective $R$-module $E$ such that $\Gid{M} = -\inf{\RHom{E}{M}}$;
  cf.~ \thmcite[3.3]{CFH-06}.  The module $E$ is a direct summand of
  an injective $S$-module $\widetilde{E}$ by \lemref{injsum}, hence
  the third step in the next sequence
  \begin{align*}
    \Gid[S]{\RHom{S}{M}}
    &\ge -\inf{\RHom[S]{\widetilde{E}}{\RHom{S}{M}}}\\
    &= -\inf{\RHom{\widetilde{E}}{M}}\\
    &\ge -\inf{\RHom{E}{M}}\\
    &=\Gid{M}.
  \end{align*}
  The first step is by \thmcite[3.3]{CFH-06}, the second one is from
  Hom-tensor adjointness, and the last one comes from the choice of
  $E$.  The opposite inequality was proved in the first paragraph of
  this proof.
\end{prf*}

The next result contains the other half of Theorem~A, and it gives a
partial answer to \cite[quest. 8.10]{SInSSW04}; see also
\prpref{SISSW2}. Its proof is similar to, but simpler than, the proof
of \thmref{01}.  Note that \thmref[]{02} has no assumptions on the
Krull dimension of $R$ or $S$.

\begin{thm}
  \label{thm:02}
  Let $\mapdef{\f}{R}{S}$ be a ring homomorphism of finite flat
  dimension.  For every homologically bounded $R$-complex $M$ there is
  an inequality
  \begin{align*}
    \Gfd{M} &\ge \Gfd[S]{\Ltpp{S}{M}}.
  \end{align*}
  If $\f$ is faithfully flat, then equality holds; in particular, the
  dimensions are simultaneously finite in this case.  \qed
\end{thm}

Equality can fail in Theorems~\pgref{thm:01} and \pgref{thm:02} if
$\f$ is not flat, even if $R$ is local and $\f$ is surjective. See
\exaref[]{nodescent} for an example.

We conclude this section with an application of \thmref{01} which, in
particular, answers \cite[quest. 8.10]{SInSSW04} for local ring
homomorphisms.

\begin{prp}
  \label{prp:SISSW2}
  Let $\mapdef{\f}{R}{S}$ be a faithfully flat ring homomorphism, and
  assume that $R$ is semi-local. For every homologically bounded
  $R$-complex $M$, there are equalities
  \begin{equation*}
    \Gfd[S]{\tpp{S}{M}}=\Gfd{\tpp{S}{M}}=\Gfd{M}.    
  \end{equation*}
\end{prp}

\begin{prf*}
  If $\Gfd{M}$ is finite, then the desired equalities hold by
  \corcite[8.9]{SInSSW04}. \thmref{02} says that $\Gfd[S]{\tpp{S}{M}}$
  and $\Gfd{M}$ are simultaneously finite. Hence, it remains to assume
  that $\Gfd{\tpp{S}{M}}$ is finite and prove that $\Gfd{M}$ is
  finite.

  The completion $\Rhat$ of $R$ (with respect to its Jacobson radical)
  has a dualizing complex.  By \thmref{02} the finiteness of
  $\Gfd{\tpp{S}{M}}$ implies that
  $\Gfd[\Rhat]{\tpp{\Rhat}{(\tp{S}{M})}}$ is finite, so the complex
  \begin{equation*}
    \tp{\Rhat}{(\tp{S}{M})}\qis \tp{\tpp{\Rhat}{M}}{S}\qis
    \tp[\Rhat]{\tpp{\Rhat}{M}}{\tpp{\Rhat}{S}}
  \end{equation*}
  is in the Auslander class $\A[\Rhat]$ by \eqref{A}. As $S$ is
  faithfully flat over $R$, the module $\tp{\Rhat}{S}$ is faithfully
  flat over $\Rhat$, and it follows that $\tp{\Rhat}{M}$ is in
  $\A[\Rhat]$, cf.~\cite[rmk. 4]{HHLDWh08}. Thus,
  $\Gfd[\Rhat]{\tpp{\Rhat}{M}}$ is finite by \eqref{A}, and
  \thmref{02} implies that $\Gfd{M}$ is finite.
\end{prf*}

\section{A Chouinard formula for Gorenstein injective dimension}
\label{sec:bass}

\noindent
The \emph{width} of a complex $M$ over a local ring $R$ with residue
field $k$ is defined as:
\begin{equation*}
  \wdt{M} = \inf{\Ltpp{k}{M}}.
\end{equation*}
There is an inequality $\wdt{M} \ge \inf{M}$, and equality holds if
$M$ is homologically finite, by Nakayama's lemma.  Let $N$ be another
$R$-complex; a standard application of the K\"unneth formula yields
\begin{equation}
  \label{wdt2} 
  \wdt\Ltpp{M}{N} =\wdt M+\wdt N.
\end{equation}
If $M$ is homologically bounded and of finite projective dimension,
and if $\H{N}$ is bounded above, then there is an equality
\thmcite[(4.14)(a) and (1.6)(b)]{CFF-02}:
\begin{equation}
  \label{wdt1} 
  \wdt\RHom{M}{N} =\dpt M+\wdt N -\dptR.
\end{equation}

Foxby~\cite{HBF79} defines the \emph{small support} of a complex $M$
over a noetherian ring $R$, denoted $\supp{M}$, as the set of prime
ideals $\p$ in $R$ such that the complex $M_\p$ has finite width over
$R_\p$.

\begin{lem}
  \label{lem:wdt}
  Let $J$ be a Gorenstein injective $R$-module. Then one has
  $$\dptR_\p \le \wdt[R_\p]{J_\p}$$ for every $\p$ in
  $\SpecR$, and equality holds if $\p$ is a maximal element in
  $\supp{J}$.
\end{lem}

\begin{prf*}
  Let $\p$ be given, and let $T$ be an $R_{\p}$-module of finite
  projective dimension.  Because there is an exact sequence
  $$\cdots \to A_2\to A_1\to A_0\to J_{\p}\to 0$$
  where each $A_i$ is an injective $R_{\p}$-module, a standard
  dimension shifting argument shows that $\Ext[R_\p]{i}{T}{J_\p}=0$
  for all $i > 0$.
  
  Set $d=\dptR_\p$, and choose a maximal $R_\p$-regular sequence
  $\x$. Because The $R_{\p}$-module $R_\p/(\x)$ has finite projective
  dimension, the previous paragraph provides the first inequality in
  the next display
  \begin{equation*}
    0 \le \inf{\RHom[R_\p]{R_\p/(\x)}{J_\p}} 
    \le  \wdt[R_\p]{\RHom[R_\p]{R_\p/(\x)}{J_\p}} 
    = \wdt[R_\p]{J_\p} -d
  \end{equation*}
  where the equality follows from \eqref{wdt1}. This proves the
  desired inequality.

  Let $\p$ be maximal in $\supp{J}$, and let $I$ be the minimal
  injective resolution of $J$. For prime ideals $\q$ that strictly
  contain $\p$, the indecomposable module $\E{R/\q}$ is not a direct
  summand of any module $I_j$ in $I$; see \cite[rmk.~2.9]{HBF79}.  It
  follows that $I_j\cong (I_j)_\p\oplus I_j'$ where $I_j'$ is a direct
  sum of injective hulls of the form $\E{R/\q}$ such that
  $\p\not\subseteq \q$. Recall that for each such $\q$ we have
  $\Hom{\E{R/\p}}{\E{R/\q}}=0$, and so $\Hom{\E{R/\p}}{I_j'}=0$.  In
  conclusion, there are isomorphisms
  $$\Hom{\E{R/\p}}{I_j} \is \Hom{\E{R/\p}}{(I_j)_\p\oplus I_j'}\is \Hom{\E{R/\p}}{(I_j)_\p}.$$
  This explains the last isomorphism below; the first one is
  Hom-tensor adjointness
  \begin{equation*}
    \Hom[R_\p]{\E[R_\p]{R_{\p}/\p R_{\p}}}{I_\p} \is \Hom{\E{R/\p}}{I_\p} \is \Hom{\E{R/\p}}{I}.
  \end{equation*}
  It follows that the modules $\Ext[R_\p]{i}{\E[R_\p]{R_{\p}/\p
      R_{\p}}}{J_\p}$ vanish for $i>0$.
  
  Set $S=R_{\p}$; it is a local ring with depth $d$, maximal ideal
  $\n:=\p R_\p$ and residue field $l:=R_{\p}/\p R_{\p}$. The
  $S$-module $B:=J_\p$ has minimal injective resolution $H:=I_\p$. One
  has $\n\in\supp[S]{B}$ and
  \begin{equation*}
    \tag{1}
    \begin{aligned}
      \Ext[S]{i}{T}{B}&=0 \ \text{for all $i>0$ and every $S$-module
        $T$
        with $\pd[S]{T}$ finite}\\
      \smash{\Ext[S]{i}{E}{B}} &=0 \ \text{for all $i>0$, where $E$ is
        the injective envelope of $l$.}
    \end{aligned}
  \end{equation*}
  To prove the desired equality $\wdt{B}=d$, we adapt the proof of
  \corcite[6.5]{LWCHHl09a}. Let $K$ denote the Koszul complex on a
  system of generators for $\n$, and note that $\tp[S]{K}{E}$ and
  $\Hom[S]{K}{E}$ are isomorphic up to a shift.  The total homology
  module $\H{\Hom[S]{K}{E}}$ has finite length. In particular
  $\tp[S]{K}{E}$ is homologically finite. Fix a resolution by finitely
  generated free $S$-modules
  \begin{equation*}
    \tag{2}
    L \qra \tp[S]{K}{E}.
  \end{equation*}
  Then there are (quasi)isomorphisms:
  \begin{equation*}
    \tag{3}
    \tp[S]{K}{\Ltpp[S]{E}{\Hom[S]{E}{B}}} 
    \qis \tp[S]{L}{\Hom[S]{E}{B}}
    \is \Hom[S]{\Hom[S]{L}{E}}{B}
  \end{equation*}
  the last one is Hom-evaluation~\lemcite[4.4]{LLAHBF91}.  The
  resolution $(2)$ induces a quasiisomorphism $\alpha$ from the
  complex $\Hom[S]{\tp[S]{K}{E}}{E} \is \Hom[S]{K}{\Shat}$ to
  $\Hom[S]{L}{E}$. The mapping cone $C$ of $\alpha$ is a bounded
  complex of direct sums of $\Shat$ and $E$. By $(1)$ the modules
  $C_j$ are Ext-orthogonal to $B$, that is, we have
  $\Ext{i}{C_j}{B}=0$ for all $i\geq 1$ and all $j$. Hence, an
  application of $\Hom[S]{-}{B}$ yields a quasiisomorphism
  \begin{equation*}
    \tag{4}
    \Hom[S]{\Hom[S]{L}{E}}{B} \xra[\qis]{\Hom[]{\alpha}{B}}
    \Hom[S]{\Hom[S]{K}{\Shat}}{B}.
  \end{equation*}
  The modules in the complex $\Hom[S]{K}{\Shat}$ are Ext-orthogonal to
  the modules in the mapping cone of the injective resolution $B \qra
  H$. Therefore, one has
  \begin{equation*}
    \tag{5}
    \Hom[S]{\Hom[S]{K}{\Shat}}{B}  \qis   \Hom[S]{\Hom[S]{K}{\Shat}}{H}
  \end{equation*}
  see \lemcite[2.4]{CFH-06}.  Now piece together $(3)$--$(5)$, and use
  Hom-evaluation to obtain
  \begin{equation*}
    \tag{6}
    \tp[S]{K}{\Ltpp[S]{E}{\Hom{E}{B}}} \qis \tp[S]{K}{\RHom[S]{\Shat}{B}}.
  \end{equation*}
  By the width sensitivity of $K$, see \cite[(4.2) and
  (4.11)]{CFF-02}, the complexes $\RHom[S]{\Shat}{B}$ and
  $\Ltp[S]{E}{\Hom[S]{E}{B}}$ have the same width. From \eqref{wdt2}
  and \eqref{wdt1} one has
  \begin{equation*}
    \tag{7}
    \wdt[S]{E} + \wdt[S]{\Hom[S]{E}{B}} = \wdt[S]{B}.
  \end{equation*}
  The maximal ideal $\n$ is in $\supp[S]{B}$, so $\wdt[S]{B}$ is
  finite. It follows from $(7)$ that $\wdt[S]{\Hom[S]{E}{B}}$ is
  finite; in particular, $\Hom[S]{E}{B}$ is non-zero.  As every
  element in $E$ is annihilated by a power of the maximal ideal $\n$,
  it follows that $\n\Hom[S]{E}{B} \ne \Hom[S]{E}{B}$.  (Indeed, if
  $\Hom[S]{E}{B}=\n\Hom[S]{E}{B}$, then
  $\Hom[S]{E}{B}=\n^t\Hom[S]{E}{B}$ for each $t\ge 1$. Since
  $\Hom[S]{E}{B}\ne 0$, there are elements $\psi\in\Hom[S]{E}{B}$ and
  $e\in E$ such that $\psi(e)\neq 0$. Also, there is an integer $t\ge
  1$ such that $\n^te=0$. The condition $\psi\in\n^t\Hom[S]{E}{B}$
  then implies $\psi(e)=0$, a contradiction.) Thus, one has
  $\wdt[S]{\Hom[S]{E}{B}}=0$, and the desired equality follows as
  $\wdt[S]{E} = d$ by \prpcite[(4.8)]{CFF-02}.
\end{prf*}

The next result contains Theorem~C from the introduction.

\begin{thm}
  \label{thm:Chouinard}
  For every $R$-complex $M$ of finite Gorenstein injective dimension
  there is an equality
  \begin{displaymath}
    \Gid{M}=\sup\{\dpt[]{R_\p}-\wdt[R_\p]{M_\p} \,|\, \p \in \SpecR\}.
  \end{displaymath}
\end{thm}

\begin{prf*}
  If $\H{M}=0$ then the equality holds for trivial reasons. Assume
  $\H{M}\ne 0$; without loss of generality, assume also that $M_0 \ne
  0$ and $M_i=0$ for all $i > 0$. Set $g=\Gid{M}$, and notice that
  $g\ge 0$.  If $g=0$, then $M$ is a Gorenstein injective module, and
  the desired equality follows immediately from \lemref{wdt}.

  Assume now that $g > 0$.  There is an exact triangle in $\D$
  \begin{equation*}
    J \to I \to M \to \Shift J
  \end{equation*}
  where $J$ is a Gorenstein injective module, and $I$ is a complex
  with $\id{I} = g$. This is dual to the special case $n = \inf{N} =
  0$ of \thmcite[3.1]{LWCSIn07}. By the Chouinard formula for
  injective dimension~\thmcite[2.10]{SYs98b}, there is a prime ideal
  $\p$ such that $\wdt[R_\p]{I_\p} = \dptR_\p - g$. By \lemref{wdt}
  one has
  \begin{equation*}
    \wdt[R_\p]{J_\p} \ge \dptR_\p > \wdt[R_\p]{I_\p}
  \end{equation*}
  so from the exact sequence of homology modules
  \begin{equation*}
    \tag{1}
    \cdots \to \H[i+1]{\Ltp{M}{R_{\p}/\p R_{\p}}} \to  \H[i]{\Ltp{J}{R_{\p}/\p R_{\p}}} \to
    \H[i]{\Ltp{I}{R_{\p}/\p R_{\p}}} \to 
    \cdots
  \end{equation*}
  one gets the equality $\wdt[R_\p]{M_\p} = \wdt[R_\p]{I_\p}$. This
  proves the inequality "$\le$".
  
  For the opposite inequality, let a prime $\q$ be given. If
  $\wdt[R_\q]{M_\q} \ge \dptR_\q$, then $g > \dptR_\q -
  \wdt[R_\q]{M_\q}$; so assume $\wdt[R_\q]{M_\q} < \dptR_\q$. Again
  $(1)$ yields $\wdt[R_\q]{M_\q} = \wdt[R_\q]{I_\q}$, as one has
  $\wdt[R_\q]{J_\q} \ge \dptR_\q$ by \lemref{wdt}. Now the inequality
  $g \ge \dptR_\q - \wdt[R_\q]{M_\q}$ follows from the Chouinard
  formula for injective dimension.
\end{prf*}

For modules, the Bass formula below is proved in \corcite[2.5]{KTY-}.
Our argument is similar; the key tools are \thmcite[3.6]{HBFAJF07} and
\thmref{Chouinard}.

\begin{cor}
  \label{cor:Bass}
  Let $R$ be local, and let $M$ be a homologically finite
  $R$-complex. If $M$ has finite Gorenstein injective dimension, then
  there is an equality $$\Gid{M} = \dptR - \inf{M}.$$
\end{cor}

\begin{prf*}
  Let $\p$ be a prime ideal in $R$, and choose a prime ideal $\q$ in
  $\Rhat$ minimal over $\p\Rhat$. The map $R_\p \to \Rhat_\q$ is local
  and flat with artinian closed fiber $\Rhat_\q/\p\Rhat_\q$. Hence one
  has $\dptR_\p - \inf{M_\p} = \dpt[]{\Rhat_\q} -
  \inf{\tpp{\Rhat}{M}_\q}$, and from \thmref{Chouinard} follows the
  inequality $$\Gid{M} \le \Gid[\Rhat]{\tpp{\Rhat}{M}}.$$ By
  \thmcite[3.6]{HBFAJF07} the complex $\tp{\Rhat}{M}$ has finite
  Gorenstein injective dimension over $\Rhat$. Since $\Rhat$ has a
  dualizing complex, one has $$\Gid[\Rhat]{\tpp{\Rhat}{M}} =
  \dpt[]{\Rhat} - \inf{\tpp{\Rhat}{M}}= \dptR - \inf{M}$$ by
  \thmcite[6.3]{CFH-06}. The two displays combine to establish the
  inequality ``$\le$''; the opposite one is from \thmref{Chouinard}.
\end{prf*}

\section{Module structures and vanishing of homology}
\label{sec:co}

\noindent
The first result of this section contains Theorem~B from the
introduction. Indeed, when $R$ is Gorenstein, every $R$-module has
finite Gorenstein injective dimension, cf.~\thmcite[(6.2.7)]{lnm}, and
\pgref{thm:FSWW} applies to the natural map $R\to\Rhat$.

\begin{thm}
  \label{thm:FSWW}
  Let $\mapdef{\f}{R}{S}$ be a flat local ring homomorphism such that
  the induced map $R/\m \to S/\m S$ is an isomorphism. Let $M$ be a
  finitely generated $R$-module with $\Gid{M}$ finite. If the
  $S$-module $\Ext{i}{S}{M}$ is finitely generated for every
  $i=1,\dots,\dim{M}$, then one has $\Ext{i}{S}{M}=0$ for $i\ge 1$,
  and $M$ has an $S$-module structure compatible with its $R$-module
  structure via $\f$.
\end{thm}

\begin{prf*}
  The module $\Hom{S}{M}$ is finitely generated over $R$ and hence
  over $S$; and the modules $\Ext{i}{S}{M}$ vanish for $i > \dim{M}$;
  see \cite[cor.~1.7 and proof of thm.~2.5]{FSW-08}. Thus, the
  $S$-complex $\RHom{S}{M}$ is homologically finite. In the sequence
  below the first and third equalities are from \corref{Bass}
  \begin{equation*}
    \dptR = \Gid{M} = \Gid[S]{\RHom{S}{M}} = \dpt[]{S} - \inf{\RHom{S}{M}}.
  \end{equation*}
  The second equality is from \thmref{01}.  The assumptions on $\f$
  imply that $R$ and $S$ have the same depth, whence
  $\inf{\RHom{S}{M}}=0$. This establishes the desired vanishing of
  Ext-modules, and the existence of the $S$-structure on $M$ follows
  from \thmcite[2.5]{FSW-08}.
\end{prf*}

\begin{rmk}
  \label{rmk:no}
  If $R$ is Gorenstein, then every finitely generated complete
  $R$-module (in particular, every $R$-module of finite length)
  satisfies the hypotheses of \thmref{FSWW}. See \thmcite[2.5]{FSW-08}
  or \thmcite[2.3]{MR2194381}.
\end{rmk}

The next example shows that the flatness hypothesis is necessary for
the equality in Theorems~\pgref{thm:01} and~\pgref{thm:02}.

\begin{exa}
  \label{exa:nodescent} Let $R$ be a complete Cohen-Macaulay local
  ring with a non-maximal prime ideal $\p\subset R$ such that $R_{\p}$
  is not Gorenstein. For example, the ring could be
  $R=k[\![X,Y,Z]\!]/(X^2,XY,Y^2)$ with prime ideal $\p=(X,Y)R$.

  As an $R$-module, $R_\p$ has infinite Gorenstein injective
  dimension.  Indeed, if $\Gid{R_{\p}}<\infty$, then
  $\Gid[R_{\p}]{R_{\p}}$ is finite as well by \prpcite[5.5]{CFH-06},
  and this contradicts the assumption that $R_{\p}$ is not Gorenstein;
  cf.~\thmcite[(6.3.2)]{lnm}.

  Let $\x=x_1,\ldots,x_d$ be a maximal $R$-regular sequence and set
  $S=R/(\x)$. The surjection $R\onto S$ is a homomorphism of finite
  flat dimension. The small supports of $S$ and $R_{\p}$ are disjoint,
  so \cite[lem.~2.6 and prop.~2.7]{HBF79} yields
  $\H{\Ltp{S}{R_{\p}}}=0$. The complexes $\Ltp{S}{R_{\p}}$ and
  $\RHom{S}{R_{\p}}$ are isomorphic (up to a shift) in $\D$. In
  particular, one has $$\Gfd[S]{\Ltpp{S}{R_{\p}}} = -\infty =
  \Gid[S]{\RHom{S}{R_{\p}}},$$ but $\Gfd{R_{\p}} = 0$ and
  $\Gid{R_{\p}} = \infty$.
\end{exa}

\begin{rmk}
  No finitely generated $R$-module can take the place of $R_\p$ in
  \exaref{nodescent}. Indeed, let $\mapdef{\f}{R}{S}$ be a local ring
  homomorphism, and let $M\neq 0$ be a finitely generated $R$-module.
  As $S/\m S$ and $M/\m M$ are not zero, then Nakayama's lemma yields
  $\tp{S}{M}\neq 0$, whence $\H{\Ltp{S}{M}}$ is not zero. Assume that
  $\f$ has finite flat dimension. Then \eqref{wdt1} yields
  $\H{\RHom{S}{M}} \ne 0$ because $\dpt{S}$ and $\wdt{M}$ are both
  finite.  Now \thmcite[4.8]{AFrSSW07a} yields
  $\Gfd[S]{\Ltpp{S}{M}}=\Gfd{M}$. Assuming further that $\f$ is
  module-finite, the corresponding equality
  $\Gid[S]{\RHom{S}{M}}=\Gid{M}$ is proved below.
\end{rmk}

\begin{prp}
  \label{prp:some}
  Let $\mapdef{\f}{R}{S}$ be a module-finite local ring homomorphism
  of finite flat dimension, and assume that $R$ admits a dualizing
  complex. For every homologically finite $R$-complex $M$ one then has
  \begin{equation*}
    \Gid{M} = \Gid[S]{\RHom{S}{M}}.
  \end{equation*}
\end{prp}

\begin{prf*}
  By \eqref{wdt1} one has $\inf{\RHom{S}{M}} = \dpt[]{S} +\inf{M}-
  \dptR$, so by \corref{Bass} it is sufficient to prove that $\Gid{M}$
  is finite if and only if $\Gid[S]{\RHom{S}{M}}$ is finite. The
  ``only if'' is already known from \thmref{01}, so assume that
  $\Gid[S]{\RHom{S}{M}}$ is finite.

  Let $D$ be a dualizing complex for $R$. Since the homomorphism $\f$
  is module-finite, the complex $\RHom{S}{D}$ is dualizing for $S$,
  cf.~\cite[(2.12)]{LLAHBF97}. By \corcite[6.4]{CFH-06} the complex
  $\RHom[S]{\RHom{S}{M}}{\RHom{S}{D}}$ has finite Gorenstein flat
  dimension over~$S$. Adjunction and
  Hom-evaluation~\lemcite[4.4]{LLAHBF91} yield
  \begin{align*}
    \RHom[S]{\RHom{S}{M}}{\RHom{S}{D}}
    &\qis \RHom{\RHom{S}{M}}{D} \\
    &\qis \Ltp{S}{\RHom{M}{D}}.
  \end{align*}
  It follows from \thmcite[4.8]{AFrSSW07a} that $\RHom{M}{D}$ has
  finite Gorenstein flat dimension over $R$, and therefore
  \corcite[(6.4)]{CFH-06} implies that $\Gid{M}$ is finite.
\end{prf*}

\begin{prp}
  Let $\mapdef{\f}{R}{S}$ be a module-finite local ring homomorphism
  of finite flat dimension. For every finitely generated complete
  $R$-module $M$ one has
  \begin{equation*}
    \Gid{M} = \Gid[S]{\RHom{S}{M}}.
  \end{equation*}
\end{prp}

\begin{prf*}
  Since $M$ is finitely generated and complete, it follows from
  \thmcite[2.5]{FSW-08} that $M$ is isomorphic to $\Hom{\Rhat}{M}$ and
  that one has $\Ext{i}{\Rhat}{M}=0$ for $i\geq 1$. In particular, the
  complex $\RHom{\Rhat}{M}$ is homologically finite over $\Rhat$.

  Let $\mapdef{\widehat{\f}}{\Rhat}{\Shat}$ denote the local
  homomorphism induced on completions. Since $S$ is module finite and
  has finite flat dimension over $R$, the completion $\Shat$ is module
  finite and has finite flat dimension over $\Rhat$.  \thmref{01}
  explains the first and fourth equalities in the next sequence:
  \begin{align*}
    \Gid{M}
    &=\Gid[\Rhat]{\RHom{\Rhat}{M}} \\
    &=\Gid[\Shat]{\RHom[\Rhat]{\Shat}{\RHom{\Rhat}{M}}}\\
    &=\Gid[\Shat]{\RHom[S]{\Shat}{\RHom{S}{M}}}\\
    &=\Gid[S]{\RHom{S}{M}}.
  \end{align*}
  The third equality is due to the isomorphisms
  $$\RHom[S]{\Shat}{\RHom{S}{M}}
  \qis\RHom{\Shat}{M} \qis\RHom[\Rhat]{\Shat}{\RHom{\Rhat}{M}}$$ and
  the second equality is from \prpref{some}.
\end{prf*}

\section*{Acknowledgments}

\noindent
We learned about Gorenstein dimensions and the question of
resolution-free characterizations from Hans-Bj\o rn Foxby.  It is a
pleasure to dedicate this paper~to~him.

We thank Hamid Rahmati for helpful discussions, and we thank Henrik
Holm for valuable comments on an earlier version of the paper.  We
also thank an anonymous referee for noticing an error in an earlier
version.

\bibliographystyle{amsplain}

\def\cprime{$'$}
  \newcommand{\arxiv}[2][AC]{\mbox{\href{http://arxiv.org/abs/#2}{\sf arXiv:#2
  [math.#1]}}}
  \newcommand{\oldarxiv}[2][AC]{\mbox{\href{http://arxiv.org/abs/math/#2}{\sf
  arXiv:math/#2
  [math.#1]}}}\providecommand{\MR}[1]{\mbox{\href{http://www.ams.org/mathscine%
t-getitem?mr=#1}{#1}}}
  \renewcommand{\MR}[1]{\mbox{\href{http://www.ams.org/mathscinet-getitem?mr=#%
1}{#1}}}
\providecommand{\bysame}{\leavevmode\hbox to3em{\hrulefill}\thinspace}
\providecommand{\MR}{\relax\ifhmode\unskip\space\fi MR }
\providecommand{\MRhref}[2]{%
  \href{http://www.ams.org/mathscinet-getitem?mr=#1}{#2}
}
\providecommand{\href}[2]{#2}


\end{document}